\documentclass[12pt]{amsart}
\usepackage{epsfig}
\usepackage{amssymb}

\def\Z{\mathbb{Z}}

\def\R{\mathbb{R}}
\def\0{\mathbb{R}_{\geq 0}}
\def\+{\mathbb{R}_{>0}}
\newtheorem{theorem}{Theorem}[section]

\theoremstyle{definition}

\title{Meunier Conjecture}
\date{29 August, 2012}
\begin{document}
\begin{abstract}
Fr\'ed\'eric Meunier's question about a multicolored Sperner lemma is addressed, leaving the question of connectivity for the color hypergraphs of such a multicolored simplex.  
Sperner's lemma asserts the existence of a simplex using all the colors for any vertex coloring of a subdivision of a large simplex with appropriate boundary conditions.  Meunier's questions generalizes this to the situation of having several such colorings and asserts the existence of a simplex using enough different colors from each coloring.  
\end{abstract}
\maketitle
\section{Notation}

A {\em finite simplicial complex} is a pair $(K,V)$ with $V$ a finite set of vertices and 
$K\subseteq PV=\{S|S\subseteq V\}$ a set of faces closed under inclusion, so that if $A\subseteq B\in K$ then $A\in K$.
The geometric realization of a simplicial complex $(K,V)$ is a subspace of the standard simplex in $\R^V=\{k:V\rightarrow\R\}$ denoted \[|(K,V)|=|K|=\{k:V \rightarrow \0|\sum_{v\in V}k(v)=1, k^{-1}(\+)\in K\}.\] 
If $(K,V)$ and $(K',V')$ are finite simplicial complexes then a {\em map} between them is any order preserving map $f:K\rightarrow K'$ meaning that if $A\subseteq B\in K$ then $f(A)\subseteq f(B) \in K'$.  
The map $f$ is {\em simplicial} if it is induced by a vertex map $f_0:V\rightarrow V'$ in the sense that $f(A)=\{f_0(v)|v\in A\}$.  
Geometric realization is functorial for simplicial maps so that if $f:K\rightarrow K'$ is simplicial then $|f|:|K|\rightarrow|K'|$ is defined for $k\in|K|$ by $(|f|(k))(v)=\sum_{\{u|f(u)=v\}}k(u)$.  
The map $f$ is a {\em subdivision} if there is a homeomorphism $h_f:|K|\rightarrow |K'|$ inducing $f$ in the sense that if $k\in|K|$ then $(h_f(k))^{-1}(\+)= f(k^{-1}(\+))$.  
If $f$ and $f'$ are maps from $K$ to $K'$ then write $f \leq f'$ and call $f'$ a deformation of $f$ and $f$ a specialization of $f'$  
if $f(A)\subseteq f'(A)$ for every $A\in K$.  Note that $\leq$ is a partial ordering of maps and that every subdivision has some simplicial specialization(s).  These simplicial specializations will be called {\em Sperner colorings} for the subdivision.  
If $(K,V)$ is a simplicial complex and $m\in \Z_{\geq 0}$ then $K_m=\{A\in K||A|\leq m\}$ and $(K_m,V)$ is called the $(m-1)$-skeleton of $K$.  
If $(K',V')$ is another simplicial complex then $K*K'=\{A\subseteq V\times V'|\pi_1(A)\in K,\pi_2(A)\in K'\}$ and $(K*K',V\times V')$ 
is called the join of $K$ and $K'$.   
If $V$ is a finite set $PV=(PV,V)$ is a simplex. 
Write $[a,b]=\{a, a+1, \ldots, b\}\subset\Z$.

\section{Theorem}

\begin{theorem}
If $f:K\rightarrow PN$ is a subdivision and $\{c_i\}_{i\in I}$ are Sperner colorings then for every collection $\{m_i\}_{i\in I}$ of sizes with $\sum_{i\in I}m_i=|N|-1$ there is a face of $K$ with every $|c_i(\sigma)|>m_i$ and $\cup_i c_i(A)=N$.  
\end{theorem}

This theorem is a variation on Sperner's lemma, which is the case in which $|I|=1$.  

\begin{proof}
If $x:N\times I\rightarrow \R_{\geq 0}$ write $Sx=\sum_{n\in N, i\in I}x(n,i)$,
$S_ix=\sum_{n\in N}x(n,i)$ and $S_nx=\sum_{i\in I}x(n,i)$.  

Consider two subspaces of the cone $\R_{\geq 0}^{N\times I}=\{f:N\times I\rightarrow \R_{\geq 0}\}$.  The first is a subspace of the product of simplices $|PN|^I$ and is denoted  
\[X=\cup_{i \in I}\left[|PN|^{I-\{i\}}\times|(PN)_{m_i}|^{\{i\}}\right]\]
\[=\{x:N\times I\rightarrow \R_{\geq 0}|\forall i\in I\hbox{ there is } S_ix=1\hbox{ and }\exists i\in I\hbox{ with }|x^{-1}(\+)\cap (N\times\{i\})|\leq m_i\}.\] 
The second is a subcomplex of the simplex $|P(N\times I)|=|PN|^{*I}$ and is the join of skeletons denoted 
\[Y=|*_{i\in I}(PN)_{m_i}|\] 
\[=\{y:N\times I\rightarrow \R_{\geq 0}|Sy=1\hbox{ and }\forall i\in I\hbox{ there is }|y^{-1}(\R_{>0})\cap (N\times\{i\})|\leq m_i\}.\]
For example if every $m_i=1$ and $|N|=|I|$ then both are spheres, with $X$ the boundary of an $|I|$-cube and $Y$ the boundary of the dual cross polytope.
This is an easy example to picture, but not one which arises above since $\sum m_i=|N|$ rather than $|N|-1$.

Define $J:Y\rightarrow X$ by setting \[(J(y))(n,i)=\alpha_{y}(y(n,i)+s_y(i))\] where 
\[s_y(i)={1\over |N|}\left[\max_{j\in I}\{S_jy\}- S_iy\right]\] so that all of the $S_iJ(y)$ are positive and equal and $\alpha_{y}={1\over \max_{j\in I}\{S_jy\}}$ is a normalizing constant so that $S_iJ(y)=1$.
For every $y\in Y$ there is some $i\in I$ with $s_y(i)=0$ and hence $y^{-1}(\+)\cap(N\times \{i\})=(J(y))^{-1}(\+)\cap(N\times \{i\})$ which has order at most $m_i$ so that $J(y)\in X$. 

Define $H:X\rightarrow Y$ by setting \[(H(x))(n,i)=\beta_x\max\{x(n,i)-r_x(i),0\}\] where 
\[r_x(i)=\min\{r||x^{-1}(\R_{>r})\cap (N\times\{i\})|\leq m_i\}\in\left[0,{1\over m_i}\right]\] 
and $\beta_x$ is a normalizing constant so that $SH(x)=1$.  
For every $x\in X$ there is some $i\in I$ with $r_x(i)=0$ and hence also some $n\in N$ with $(H(x))(n,i)\not=0$ so there is a choice for $\beta_x$ with $H(x)\in Y$.


Note that $J$ is an embedding, $H\circ J=\hbox{id}_Y$ and $J\circ H\sim \hbox{id}_X$ via the linear homotopy using the linear structure on $\R^{N\times I}$ with $h(x,t)=tx+(1-t)J(H(x))$.  Thus $H$ is a homotopy equivalence.  

Define $\rho:Y\rightarrow |(PN)_{|N|-1}|=\partial|PN|\cong S^{|N|-1}$ by setting $(\rho(y))(n)=S_ny$.  

If $\{c_i\}$ are from a counterexample to the theorem then 
$C=\prod_{i\in I}|c_i|:|K|\rightarrow X$.  

Since every facet $F$ of $|PN|$ has $\rho(H(C(h_f^{-1}(F))))\subseteq F$ the map $(\rho\circ H\circ C)|_{\partial|K|}$ to $\partial|PN|$ is linearly homotopic to the homeomorphism $h_f|_{\partial|K|}$ and hence a degree one map of spheres.  Since it extends to the ball $\rho\circ H\circ C:|K|\rightarrow \partial|PN|$ it is also degree zero.  
\end{proof}
\section{Examples}
{\bf Example 1.}  Consider $\Pi_{n,r}$ to be all partitions of a finite sequence of $r$ objects into $n$ intervals.  The intervals can be empty, but still have a location in the sequence.  A partition is thus equivalent to a nondecreasing map $\pi:[0,n]\rightarrow [0,r]$ with $\pi(0)=0$ and $\pi(n)=r$ by taking the $i$th interval to be $\{ t\in[1,r]|\pi(i-1)<t\leq\pi(i)\}$.  

Consider the subdivision $K_{n,r}$ of the $(n-1)$ simplex with vertices $\Pi_{n,r}$ and faces \[\{A\subseteq\Pi_{n,r}|\forall a,b\in A \exists e\in\{-1,1\}\forall i\hbox{ there is } a(i)\in\{b(i),b(i)+e\}\}.\]  Thus faces are collections of partitions so that in every pair all the intervals in one of the partitions start on either the same object as the corresponding interval in the other partition or else the next object.  The subdivision map $f_{n,r}:K_{n,r}\rightarrow P[1,n]$ is defined by $f_{n,r}(A)=\{i\in[n]|\exists a\in A\hbox{ with }a(i-1)\not=a(i)\}$ so that a collection of partitions is taken to the set of interval labels which label a nonempty interval in at least one of the partitions.  
A Sperner coloring is then a rating scheme which chooses a best interval from each partitioning so that no empty interval is ever best.    
An interval in a partition is considered close to optimal under the rating scheme if shifting some of the interval start points forward by one makes it optimal.  Several colorings arise if there are several different schemata for valuing the intervals and the theorem guarantees a decomposition with enough intervals which are close to optimal using each of several valuation schemata.  
\vskip5pt
{\bf Example 2.}  If $n=3$ then example 1 is the regular triangular decomposition of a triangle with $r+1$ vertices along each side.  If the rating scheme prefers longer intervals then the coloring is the simplicial map which assigns to each vertex of the subdivision the closest corner of the triangle.
\vskip5pt
In the situation of the theorem with $\sigma$ a face of $K$ call the multiset $\{c_i(\sigma)\}_{i\in I}$ the {\em color hypergraph} of $\sigma$.  
The theorem asserts the existence of a simplex whose color hypergraph has no isolated colors and face sizes at least as large as those of some fixed spanning hypertree.  It does not guarantee the existence of one with color hypergraph containing an isomorphic copy of that spanning hypertree, as illustrated in the following two examples.  
\vskip5pt
For the last two examples consider the above example 1 with $n=4$ intervals, three colorings, every $m_i=1$ and $r=2s+1\geq 5$.  These are thus subdivisions of a tetrahedron.  
\vskip5pt
{\bf Example 3.} Consider the colorings in which $c_i$ prefers the $i$th interval unless it is empty in which case they all prefer the $4$th interval.  If this too is empty the tie breaking scheme does not affect the example.  In this case the only solutions are some of the simplices with vertices having functions $\pi$ with $(\pi(0),\pi(1),\pi(2),\pi(3),\pi(4))$ being among: $(0,1,2,3,r)$, $(0,1,2,2,r)$, $(0,1,1,2,r)$, $(0,0,1,2,r)$, $(0,1,1,1,r)$, $(0,0,1,1,r)$, $(0,0,0,1,r)$ and $(0,0,0,0,r)$.  These are the vertices of a cube which is subdivided into six tetrahedra and the faces which are solutions are exactly those which contain the long diagonal edge with ends labeled by $(0,0,1,1,r)$ and $(0,1,1,2,r)$.  The $3$ colorings of these eight vertices are $(1,1,1,4,1,4,4,4)$, $(2,2,4,2,4,2,4,4)$ and $(3,4,3,3,4,4,3,4)$ so that the coloring hypergraph for every solution face is the tree with $4$ vertices, of which $3$ are leaves.  

\vskip5pt
{\bf Example 4.} Consider the coloring $c_1$ which prefers interval $1$ and if it is empty then interval $3$, the coloring $c_2$ which prefers interval $2$ and if it is empty then interval $4$ and the coloring $c_3$ which prefers the longest interval.  Ties can again be broken in any way.  In this case the only solutions are simplices with vertices having functions $\pi$ with $(\pi(0),\pi(1),\pi(2),\pi(3),\pi(4))$ from among: $(0,0,0,s+1,r)$, $(0,0,0,s,r)$, $(0,1,1,s+1,r)$ and $(0,0,1,s+1,r)$ and six others.  These are the vertices of a convex polytope which is subdivided into ten tetrahedra and exactly which faces are solutions will depend on the tie breaking, but the tetrahedron with the four named vertices is always among them and all will have the same color hypergraph.  The $3$ colorings of this particular simplex are then $(3,3,1,3)$, $(4,4,4,2)$ and $(3,4,3\hbox{ or }4,3\hbox{ or }4)$ so that the coloring hypergraph is a path, which is the only other tree with $4$ vertices.   
\vskip5pt
{\bf Question.} In the situation of the theorem, is there always a face satisfying the theorem and for which the color hypergraph is connected?

\bibliographystyle{amsplain}

\def\cprime{$'$} \def\cprime{$'$}
\providecommand{\bysame}{\leavevmode\hbox to3em{\hrulefill}\thinspace}
\providecommand{\MR}{\relax\ifhmode\unskip\space\fi MR }
\providecommand{\MRhref}[2]{%
  \href{http://www.ams.org/mathscinet-getitem?mr=#1}{#2}
}
\providecommand{\href}[2]{#2}

\vskip20pt
\noindent{\bf Author.} Eric Babson
\end{document}